\documentclass[a4paper,11pt]{article}
\usepackage{latexsym,amsfonts,amsmath}
\usepackage{graphics}

\def\inter{\mathop{\cap}}
\def\NN{\mathbb{N}}

\def\RR{\mathbb{R}}

\newcommand{\widebar}[1]{\overline{#1}}

\def\inter{\mathop{\cap}}
\def\liminf{\mathop{\underline{\lim}}}
\def\limsup{\mathop{\overline{\lim}}}
\def\ds{\displaystyle}

\def\nl{\mbox{} \newline }

\usepackage{natbib}

\begin{document}
\newtheorem{theorem}{Theorem}[section]
\newtheorem{proposition}[theorem]{Proposition}
\newtheorem{lemma}[theorem]{Lemma}
\newtheorem{corollary}[theorem]{Corollary}
\newtheorem{definition}[theorem]{Definition}
\newtheorem{remark}[theorem]{Remark}
\newtheorem{conjecture}[theorem]{Conjecture}
\newtheorem{assumption}[theorem]{Assumption}

\bibliographystyle{plain}

\title{The Policy Iteration Algorithm for Average Continuous Control of Piecewise Deterministic Markov Processes
}
\author{ \mbox{ }
\\
O.L.V. Costa
\thanks{This author received financial support from CNPq (Brazilian National Research
Council), grant 304866/03-2 and FAPESP (Research Council of the
State of S\~ao Paulo), grant 03/06736-7.}
\thanks{Author to whom correspondence should be sent to.}
\\ \small Departamento de Engenharia de Telecomunica\c c\~oes e
Controle \\ \small Escola Polit\'ecnica da Universidade de S\~ao
Paulo \\ \small CEP: 05508 900-S\~ao Paulo, Brazil.
\\ \small phone: 55 11 30915771; fax: 55 11 30915718.
 \\ \small
e-mail: oswaldo@lac.usp.br
\\ \and
\\ F. Dufour
\\
\small Universite Bordeaux I \\
\small IMB, Institut Math\'ematiques de Bordeaux \\
\small INRIA Bordeaux Sud Ouest, Team: CQFD \\
\small \small 351 cours de la Liberation \\
\small 33405 Talence Cedex, France \\ \small e-mail :
dufour@math.u-bordeaux1.fr }
\maketitle

\begin{abstract}
The main goal of this paper is to apply the so-called policy iteration algorithm (PIA) for the
long run average continuous control problem of piecewise deterministic Markov processes (PDMP's)
taking values in a general Borel space and with compact action space depending on the state variable.
In order to do that we first derive some important properties for a pseudo-Poisson equation associated to the problem.
In the sequence it is shown that the convergence of the PIA to a solution satisfying the optimality equation holds under some classical hypotheses and that
this optimal solution yields to an optimal control strategy for the average control problem for the continuous-time PDMP in a feedback form.
\end{abstract}
\begin{tabbing}
\small \hspace*{\parindent}  \= {\bf Keywords:}
piecewise-deterministic Markov Processes, continuous-time,
long-run average\\ cost, optimal control,
integro-differential optimality inequation, policy iteration algorithm.\\
\> {\bf AMS 2000 subject classification:} 60J25, 90C40, 93E20
\end{tabbing}

\newpage
\section{Introuction}
\label{intro}
This paper studies the policy iteration algorithm (PIA) for the average cost control problem of a class of continuous-time Markov processes, namely piecewise-deterministic Markov processes (PDMP's).
These processes have been introduced in the literature by M.H.A. Davis \cite{davis93} as a general class of stochastic models.
They are a family of Markov processes involving deterministic motion punctuated by random jumps.
The motion of the PDMP $\{X(t)\}$ depends on three local characteristics, namely the flow $\phi$, the jump rate $\lambda$ and the
transition measure $Q$, which specifies the post-jump location.
Starting from $x$ the motion of the process follows the flow $\phi(x,t)$ until the first jump time $T_1$ which occurs either
spontaneously in a Poisson-like fashion with rate $\lambda(\phi(x,t))$ or when the flow $\phi(x,t)$ hits the boundary of the state-space.
In either case the location $Z_1$ of the process at the jump time $T_1$ is selected by the transition measure $Q(\phi(x,T_1),.)$.
Starting from $Z_{1}$, we now select the next interjump time $T_{2}-T_{1}$ and postjump location $X(T_{2})=Z_{2}$ in a similar way.
This gives a piecewise deterministic trajectory for $\{X(t)\}$ with jump times $\{T_{k}\}$ and postjump locations $\{Z_{k}\}$, and which follows the flow $\phi$ between two jumps.
A suitable choice of the state space and the local characteristics $\phi$, $\lambda$, and $Q$ provide stochastic
models covering a great number of problems of operations research \cite{davis93}.


The present work is a continuation of a series of papers: \cite{average,vanishing}.
It deals with the long run average cost control problem of PDMP's taking values in a general Borel space. At each point $x$ of the state space
a control variable is chosen from a compact action set $\mathbb{U}(x)$ and is applied on the jump parameter $\lambda$ and transition measure $Q$.
The long run average cost is composed of a running cost and a boundary cost (which is added each time the PDMP touches the boundary).
In this context, we follow the idea developed in \cite{average,vanishing} consisting of writing the
optimality equation for the long run average cost control problem of the PDMP $\{X(t)\}$
in terms of a discrete-time optimality equation related to the embedded Markov chain given by the post-jump location of the process $\{X(t)\}$.
As pointed out in \cite{average}, this discrete-time optimality equation is different from those classical ones encountered within the context of discrete-time Markov decision processes.
The two main reasons for doing that is to use the powerful tools developed in the discrete-time framework
(see for example the references \cite{bertsekas78,dynkin79,hernandez96,hernandez99})
and to avoid working with the infinitesimal generator associated to a PDMP, which in most cases has its domain of definition
difficult to be characterized.

The PIA has received considerable attention in the literature and consists of three steps: initialization, policy evaluation, which is related to the Poisson equation (PE) associated to the transition law
defining the Markov decision process, and policy improvement. Without attempting to present here an exhaustive panorama of the literature for the PIA, we can mention
the surveys \cite{arapostathis93,borkar91,hernandez91,hernandez99,puterman94} and the references therein and more specifically
the references \cite{PIA-97,meyn97} that analyze in details the PIA for general Markov decision processes and provide conditions which guarantee its converge.

The paper is organized as follows. We shall formulate in section \ref{formulation} the control problem while in section \ref{AssDef} some of the main assumptions are presented.
In our context, the policy evaluation step is connected to a kind of PE which we call a pseudo-Poisson equation.
This equation is clearly different from a classical PE encountered in the literature of the discrete-time Markov control processes, see Remark \ref{difference}.
However, although different, we can show in section \ref{pseudo} that this pseudo-Poisson equation still has the
good properties that we might expect to satisfy in order to guarantee the convergence of the policy iteration algorithm.
These results are not straightforward to obtain due to the specific structure of this discrete-time optimality equation.
Finally in section \ref{PIA}, the PIA is studied in details.
It is first shown that the convergence of the PIA to a solution satisfying the optimality equation holds under some classical hypotheses. In the sequence it is shown that
this optimal solution yields to an
optimal control strategy for the average control problem for the continuous-time PDMP in a feedback form.

\section{Definitions and problem formulation}
\label{formulation}
\subsection{Presentation of the control problem}
\label{pre}
In this section we present some standard notation and some basic definitions related to the motion of a PDMP $\{X(t)\}$,
and the control problems we will consider throughout the paper.
For further details and properties the reader is referred to \cite{davis93}. The following notation will be used in this paper:
$\NN$ denotes the set of natural numbers, $\RR$ the set of real numbers, $\RR_+$ the set of positive real numbers and $\RR^d$ the $d$-dimensional euclidian space. We write
$\eta$ as the Lebesgue measure on $\RR$. For $X$ a metric space $\mathcal{B}(X)$ represents the $\sigma$-algebra generated by the open sets of $X$. $\mathcal{M}(X)$ (respectively, $\mathcal{P}(X)$)
denotes the set of all finite (respectively probability) measures on $(X,\mathcal{B}(X))$. Let $X$ and $Y$ be metric spaces.
The set of all Borel measurable (respectively bounded) functions from $X$ into $Y$ is denoted by $\mathbb{M}(X;Y)$ (respectively $\mathbb{B}(X;Y)$).
Moreover, for notational simplicity $\mathbb{M}(X)$ (respectively $\mathbb{B}(X)$, $\mathbb{M}(X)^{+}$, $\mathbb{B}(X)^{+}$) denotes $\mathbb{M}(X;\RR)$
(respectively $\mathbb{B}(X;\RR)$, $\mathbb{M}(X;\RR_{+})$, $\mathbb{B}(X;\RR_{+})$).
For $g\in \mathbb{M}(X)$ with $g(x)>0$ for all $x\in X$, $\mathbb{B}_{g}(X)$ is the set of functions $v\in \mathbb{M}(X)$ such that $\ds ||v(x)||_{g}=\sup_{x\in X} \frac{|v(x)|}{g(x)}< +\infty$.
$\mathbb{C}(X)$ denotes the set of continuous functions from $X$ to $\RR$.
For $h\in \mathbb{M}(E)$, $h^+$ (respectively $h^{-}$) denotes the positive (respectively, negtive) part of $h$.

\noindent
Let $E$ be an open subset of $\RR^n$, $\partial E$ its boundary, and $\widebar{E}$ its closure. A controlled PDMP is determined by its local
characteristics $(\phi,\lambda,Q)$, as presented in the sequel.
The flow $\phi(x,t)$ is a function $\phi: \: \mathbb{R}^{n}\times \RR_{+} \longrightarrow \mathbb{R}^{n}$ continuous in $(x,t)$ and such that $\phi(x,t+s) = \phi(\phi(x,t),s).$
For each $x\in E$ the time the flow takes to reach the boundary starting from $x$ is defined as
$t_{*}(x)\doteq \inf \{t>0:\phi(x,t)\in \partial E \}$. For $x\in E$ such that $t_{*}(x)=\infty$ (that is, the flow starting from $x$ never touches the boundary), we set
$\phi(x,t_{*}(x))=\Delta$, where $\Delta$ is a fixed point in $\partial E$.
 We define the following space of functions absolutely continuous along the flow with limit towards the boundary:
\begin{align*}
\mathbb{M}^{ac}(E) & =  \bigl\{ g\in\mathbb{M}(E) \: : \: g(\phi(x,t)): [0,t_{*}(x)) \mapsto \RR \text{ is absolutely continuous for each } x\in E\\
& \text{ and whenever } t_{*}(x)< \infty \text{ the limit }\lim_{t\rightarrow t_{*}(x)} g(\phi(x,t)) \text{ exists} \bigr\}.
\end{align*}
For $g\in \mathbb{M}^{ac}(E)$ and $z\in \partial E$ for which there exists $x\in E$ such that $z=\phi(x,t_{*}(x))$ where $t_{*}(x)< \infty$ we define $\ds g(z) = \lim_{t\rightarrow t_{*}(x)} g(\phi(x,t))$
(note that the limit exists by assumption). As shown in Lemma 2 in \cite{ECC07}, for $g\in\mathbb{M}^{ac}(E)$ there exists a function $\mathcal{X}g \in \mathbb{M}(E)$ such that for all $x\in E$ and $t\in [0,t_{*}(x))$
$g(\phi(x,t))-g(x)  =  \int_{0}^{t} \mathcal{X}g(\phi(x,s)) ds $.

\bigskip

\noindent
The local characteristics $\lambda$ and $Q$ depend on a control action $u\in \mathbb{U}$ where $\mathbb{U}$ is a compact metric space
(there is no loss of generality in assuming this property for $\mathbb{U}$, see Remark 2.8 in \cite{average}), in the following way:
$\lambda \in \mathbb{M}(\widebar{E}\times\mathbb{U})^{+}$ and 
$Q$ is a stochastic kernel on $E$ given $\widebar{E}\times \mathbb{U}$.
For each $x\in \widebar{E}$ we define the subsets $\mathbb{U}(x)$ of $\mathbb{U}$ as the set of feasible control actions that can be taken when the state process is $x\in \widebar{E}$,
that is, the control action that will be applied to $\lambda$ and $Q$ must belong to $\mathbb{U}(x)$.
The following assumptions, based on the standard theory of Markov decision processes (see for example \cite{hernandez96}), will be made throughout the paper:
\begin{assumption}
\label{Hyp1a} For all $x\in \widebar{E}$, $\mathbb{U}(x)$ is a compact subspace of $\mathbb{U}$.
\end{assumption}
\begin{assumption}
\label{Hyp2a} The set $K=\left\{(x,a): x\in \widebar{E}, a \in \mathbb{U}(x) \right\}$ is a Borel subset of $\widebar{E}\times \mathbb{U}$.
\end{assumption}

\noindent We present next the definition of an admissible control strategy and the associated motion of the controlled process.
A control policy $U$ is a pair of functions $(u,u_{\partial}) \in \mathbb{M}(\NN \times E\times \RR_{+};\mathbb{U}) \times \mathbb{M}(\NN \times E;\mathbb{U})$
satisfying $u(n,x,t)\in \mathbb{U}(\phi(x,t))$, and $u_{\partial}(n,x)\in \mathbb{U}(\phi(x,t_{*}(x)))$ for all $(n,x,t)\in \NN \times E\times \RR_{+}$.
The class of admissible control strategies will be denoted by $\mathcal{U}$.
Consider the state space $\widehat{E}=E\times E \times \RR_{+} \times \NN$.
For a control policy $U=(u,u_{\partial})$ let us introduce the following parameters for $\hat{x}=(x,z,s,n)\in \widehat{E}$:
the flow $\widehat{\phi}(\hat{x},t) = (\phi(x,t), z , s+t , n)$,
the jump rate $\widehat{\lambda}^{U}(\hat{x})=\lambda(x,u(n,z,s))$, and
the transition measure
\begin{eqnarray*}
\widehat{Q}^{U}(\hat{x},A\times B\times \{0\}\times \{n+1\}) =
\begin{cases}
Q(x,u(n,z,s)); A\inter B) & \text{ if } x\in E,\\
Q(x,u_{\partial}(n,z);A\inter B) & \text{ if } x\in \partial E,
\end{cases}
\end{eqnarray*}
for $A$ and $B$ in $\mathcal{B}(E)$. From \cite[section 25]{davis93}, it can be shown that for any control strategy $U=(u,u_{\partial})\in \mathcal{U}$
there exists a filtered probability space $(\Omega,\mathcal{F},\{ \mathcal{F}_{t} \}, \{ P^{U}_{\hat{x}} \}_{\hat{x}\in \widehat{E}})$
such that the piecewise deterministic Markov process $\{\widehat{X}^{U}(t)\}$ with local characteristics $(\widehat{\phi},\widehat{\lambda}^{U},\widehat{Q}^{U})$ may be constructed as follows.
For notational simplicity the probability $P^{U}_{\hat{x}_{0}}$ will be denoted by $P^{U}_{(x,k)}$ for $\hat{x}_{0}=(x,x,0,k)\in \widehat{E}$.
Take a random variable $T_1$ such that
\begin{equation*}
P^{U}_{(x,k)}(T_1>t) \doteq
\begin{cases}
e^{-\Lambda^{U}(x,k,t)} & \text{for } t<t_{*}(x)\\
0 & \text{for } t\geq t_{*}(x)
\end{cases}
\end{equation*}
where for $x\in E$ and $t\in [0,t_*(x)[$, $\Lambda^{U}(x,k,t) \doteq \int_0^t\lambda(\phi(x,s),u(k,x,s))ds.$
If $T_1$ is equal to infinity, then for $t\in \RR_+$, $\widehat{X}^{U}(t)= \bigl(\phi(x,t),x,t,k\bigr)$.
Otherwise select independently an $\widehat{E}$-valued random variable (labelled $\widehat{X}^{U}_{1}$) having distribution
\begin{align*}
P^{U}_{(x,k)}(\widehat{X}^{U}_{1} \in A \times B \times \{0\}
\times \{k+1\} |\sigma\{T_1\}) =
\begin{cases}
Q(\phi(x,T_{1}),u(k,x,T_{1})); A\inter B) & \text{ if } \phi(x,T_{1})\in E, \\
Q(\phi(x,T_{1}),u_{\partial}(k,x); A\inter B) & \text{ if }
\phi(x,T_{1}) \in \partial E.
\end{cases}
\end{align*}
The trajectory of $\{\widehat{X}^{U}(t)\}$ starting from $(x,x,0,k)$, for $t\leq T_1$ , is given by
\begin{equation*}
\widehat{X}^{U}(t) \doteq
\begin{cases}
\bigl(\phi(x,t),x,t,k\bigr) &\text{for } t<T_1, \\
\widehat{X}^{U}_1 &\text{for } t=T_1.
\end{cases}
\end{equation*}
Starting from $\widehat{X}^{U}(T_1)=\widehat{X}^{U}_1$, we now select the next inter-jump time $T_2-T_1$ and post-jump location $\widehat{X}^{U}(T_2)=\widehat{X}^{U}_2$ in a similar way.
Let us define the components of the PDMP $\{\widehat{X}^{U}(t)\}$ by
\begin{eqnarray}
\widehat{X}^{U}(t)=\bigl(X(t),Z(t),\tau(t),N(t)\bigr).
\label{defXU}
\end{eqnarray}
For notational convenience, we have omitted to write explicitly the dependence of $U$ on the components: $X(t)$, $Z(t)$, $\tau(t)$ and $N(t)$.
From the previous construction, it is easy to see that $X(t)$ corresponds to the trajectory of the
system, $Z(t)$ is the value of $X(t)$ at the last jump time before $t$, $\tau(t)$ is the time elapsed from the last jump up to time $t$,
and $N(t)$ is the number of jumps of the process  $\{X(t)\}$ up to time $t$. As in Davis \cite{davis93}, we consider the following assumption to avoid any accumulation point of the jump times:
\begin{assumption}
\label{Hypjump} For any $x\in E$, $U\in \mathcal{U}$, and $t\geq 0$, we have $\ds E^{U}_{(x,0)}\Biggl[
\sum_{i=1}^{\infty} I_{\{T_{i}\leq t\}} \Biggr] < \infty$.
\end{assumption}
\begin{remark}
\label{RemJump}
In particular, a consequence of Assumption \ref{Hypjump} is that $T_{m}\rightarrow \infty$ as $m \rightarrow \infty$ $P^U_{(x,0)}$ for all $x\in E$, $U\in \mathcal{U}$.
\end{remark}

\noindent
The costs of our control problem will contain two terms, a running cost $f$ and a boundary cost $r$, satisfying the following properties:
\begin{assumption}
\label{Hyp6a} $f\in \mathbb{M}(\widebar{E}\times\mathbb{U})^{+}$, and $r\in \mathbb{M}(\partial E\times\mathbb{U})^{+}$.
\end{assumption}
Define for $\alpha \geq 0$, $t\in \mathbb{R}_+$, and $U\in \mathcal{U}$,
\begin{align*}
\mathbf{J}^\alpha(U,t) = \int_{0}^{t} e^{-\alpha s}f\bigl(X(s), & u(N(s),Z(s),\tau(s)) \bigr) ds  + \int_{0}^{t}e^{-\alpha s} r\bigl(X(s-),u_{\partial}(N(s-),Z(s-)) \bigr)
dp^{*}(s),
\end{align*}
where $\ds p^{*}(t) = \sum_{i=1}^{\infty} I_{\{T_{i}\leq t\}} I_{\{ X(T_{i}-) \in \partial E\}}$ counts the number of times the
process hits the boundary up to time $t$ and, for notational simplicity, set $\mathbf{J}(U,t)=\mathbf{J}^0(U,t)$.
The long-run average cost we want to minimize over $\mathcal{U}$ is given by:
$\ds \mathcal{A}(U,x)=\limsup_{t\rightarrow +\infty} \frac{1}{t} E^{U}_{(x,0)} [\mathbf{J}(U,t)]$.
We need the following assumption, to avoid infinite costs for the discounted case, see \cite{average}.
\begin{assumption}
\label{Hypdis} For all $\alpha>0$ and all $x\in E$, $\ds \inf_{U\in\mathcal{U}} E^{U}_{(x,0)}[\mathbf{J}^\alpha(U,\infty)]<\infty$.
\end{assumption}

\subsection{Discrete-time relaxed and ordinary controls}
\label{conrex1}
We present in this sub-section the set of discrete-time relaxed and ordinary controls.
Consider $\mathbb{C}(\mathbb{U})$ equipped with the topology of uniform convergence and $\mathcal{M}(\mathbb{U})$ equipped with the
weak$^{*}$ topology $\sigma(\mathcal{M}(\mathbb{U}),\mathbb{C}(\mathbb{U}))$.
For $x\in E$, define $\mathcal{P}_{x}\bigl(\mathbb{U}\bigr)$ as the set of measures $\mu\in\mathcal{P}(\mathbb{U})$ satisfying $\mu(\mathbb{U}(\phi(x,t_{*}(x))))=1$.
$\mathcal{P}(\mathbb{U})$ and $\mathcal{P}_{x}(\mathbb{U})$ for $x\in E$ are subsets of $\mathcal{M}(\mathbb{U})$ and are equipped with the relative topology.

Let $\mathcal{V}^{r}$ (respectively $\mathcal{V}^{r}(x)$ for $x\in E$) be the set of all $\eta$-measurable functions $\mu$ defined on $\RR_{+}$ with value in $\mathcal{P}(\mathbb{U})$ such that
$\mu(t,\mathbb{U})=1$ $\eta$-a.e. (respectively $\mu(t,\mathbb{U}(\phi(x,t)))=1$ $\eta$-a.e.).
It can be shown (see sub-section 3.1 in \cite{average}) that $\mathcal{V}^{r}(x)$ is a compact set of the metric space $\mathcal{V}^{r}$: a sequence $\bigl(\mu_{n}\bigr)_{n\in \NN}$ in $\mathcal{V}^{r}(x)$ converges to $\mu$ if and only if for all $g\in L^{1}(\RR_{+};\mathbb{C}(\mathbb{U}))$
\begin{eqnarray*}
\lim_{n\rightarrow \infty} \int_{\RR_{+}} \int_{\mathbb{U}(\phi(x,t))}  g(t,u) \mu_{n}(t,du) dt = \int_{\RR_{+}} \int_{\mathbb{U}(\phi(x,t))}  g(t,u) \mu(t,du) dt.
\end{eqnarray*}
The sets of relaxed controls can be defined as follows:
$\mathbb{V}^{r}(x)  = \mathcal{V}^r(x) \times \mathcal{P}_{x}\bigl(\mathbb{U}\bigr)$, for $x\in E$ 
and $\mathbb{V}^{r}  = \mathcal{V}^r \times \mathcal{P}\bigl(\mathbb{U}\bigr)$.
The set of ordinary controls, denoted by $\mathbb{V}$ (respectively $\mathbb{V}(x)$ for $x\in E$), is defined as above
except that it is composed of deterministic functions instead of probability measures. More specifically we have
$\mathcal{V}(x)  =  \bigl\{ \nu \in \mathbb{M}(\RR_{+}, \mathbb{U}) : (\forall t\in \RR_{+}), \nu(t) \in \mathbb{U}(\phi(x,t)) \bigr\}$,
$\mathbb{V}(x)  =  \mathcal{V}(x) \times \mathbb{U}(\phi(x,t_{*}(x))),$
$\mathbb{V}  =  \mathbb{M}(\RR_{+}, \mathbb{U}) \times
\mathbb{U}$.
Consequently, the set of ordinary controls is a subset of the set of relaxed controls $\mathbb{V}^{r}$ (respectively
$\mathbb{V}^{r}(x)$ for $x\in E$) by identifying any control action $u\in \mathbb{U}$ with the Dirac measure concentrated on
$u$. Thus we can write that $\mathbb{V}\subset \mathbb{V}^{r}$ (respectively $\mathbb{V}(x)\subset \mathbb{V}^{r}(x)$ for $x\in
E$) and from now on we will consider that $\mathbb{V}$ (respectively $\mathbb{V}(x)$ for $x\in E$) will be endowed with
the topology generated by $\mathbb{V}^{r}$.
The necessity to introduce the class of relaxed control $\mathbb{V}^{r}$ is justified by the fact that in general there
does not exist a topology for which $\mathbb{V}$ and $\mathbb{V}(x)$ are compact sets.

\noindent
As in \cite{hernandez96}, page 14, we need that the set of feasible state/relaxed-control pairs is a measurable subset of $\mathcal{B}(E)\times \mathcal{B}(\mathbb{V}^{r})$, that is, we need the following assumption.
\begin{assumption}
\label{Mesurability}
$\mathcal{K} \doteq \bigl\{ (x,\Theta) : \Theta \in \mathbb{V}^{r}(x), x\in E \bigr\} \in \mathcal{B}(E)\times \mathcal{B}(\mathbb{V}^{r}).$
\end{assumption}
A sufficient condition is presented in \cite[Proposition 3.3]{average} to ensure that Assumption \ref{Mesurability} holds.

\subsection{Discrete-time operators and measurability properties}
In this sub-section we present some important operators associated to the optimality equation of the discrete-time problem.
We consider the following notation $\ds w(x,\mu)  \doteq  \int_{\mathbb{U}} w(x,u)  \mu (du)$ and
$\ds Qh(x,\mu) \doteq \int_{\mathbb{U}} \int_{E} h(z) Q(x,u;dz) \mu (du)$, and
$\ds \lambda Qh(x,\mu) \doteq  \int_{\mathbb{U}} \lambda(x,u) \int_{E} h(z) Q(x,u;dz)  \mu (du)$
for $x\in \widebar{E}$, $\mu\in\mathcal{P}\bigl(\mathbb{U}\bigr)$, $h\in \mathbb{M}(E)^{+}$ and $w\in \mathbb{M}(\widebar{E}\times \mathbb{U})^{+}$.

The following operators will be associated to the optimality equations of the discrete-time problems that will be
presented in the next sections. For $\Theta=\bigl(\mu,\mu_{\partial}\bigr)\in \mathbb{V}^{r}$, $(x,A)\in E\times \mathcal{B}(E)$, $\alpha \in \RR$, according to Lemma 2 in \cite[Appendix 5]{dynkin79}
define
\begin{eqnarray}
\Lambda^{\mu}(x,t)  & \doteq & \int_{0}^{t} \lambda(\phi(y,s),\mu (s)) ds \nonumber \\
\label{DefGr}
G_{\alpha}(x,\Theta;A) & \doteq & \int_0^{t_{*}(x)}e^{-\alpha s - \Lambda^{\mu}(x,s)}\lambda QI_{A}(\phi(x,s),\mu(s)) ds\nonumber \\
& &+ e^{-\alpha t_{*}(x) -\Lambda^{\mu}(x,t_{*}(x))} Q(\phi(x,t_{*}(x)),\mu_{\partial};A).
\end{eqnarray}
For $h\in \mathbb{M}(E)^{+}$, we define $G_{\alpha}h(x,\Theta)\doteq\ds \int_{E} h(y) G_{\alpha}(x,\Theta;dy)$.
For $x\in E$, $\Theta=\bigl(\mu,\mu_{\partial}\bigr)\in \mathbb{V}^{r}$, $v\in \mathbb{M}(E\times \mathbb{U})^{+}$, $w\in \mathbb{M}(\partial E\times \mathbb{U})^{+}$,
$\alpha \in \RR$, introduce
\begin{eqnarray}
\label{DefLr}
L_{\alpha}v(x,\Theta) & \doteq  & \int_0^{t_{*}(x)}e^{-\alpha s-\Lambda^{\mu}(x,s)} v(\phi(x,s),\mu(s)) ds, \\
\label{DefHr} H_{\alpha}w(x,\Theta) & \doteq & e^{-\alpha t_{*}(x)-\Lambda^{\mu}(x,t_{*}(x))} w(\phi(x,t_{*}(x)),\mu_{\partial}).
\end{eqnarray}

For $h\in \mathbb{M}(E)$ (respectively, $v\in \mathbb{M}(E\times \mathbb{U})$),  $G_{\alpha}h(x,\Theta)=G_{\alpha}h^{+}(x,\Theta)-G_{\alpha}h^{-}(x,\Theta)$
(respectively, $L_{\alpha}v(x,\Theta)=L_{\alpha}v^{+}(x,\Theta)-L_{\alpha}v^{-}(x,\Theta)$) provided the difference has a meaning.
It will be useful in the sequel to define the function $\mathcal{L}_{\alpha}(x,\Theta)$  as follows:
$\mathcal{L}_{\alpha}(x,\Theta)  \doteq  L_{\alpha}I_{E\times \mathbb{U}}(x,\Theta)$. In particular for $\alpha =0$ we write for simplicity $G_0=G$, $L_0=L$, $H_0=H$, $\mathcal{L}_0=\mathcal{L}$.
Measurability properties of the operators $G_{\alpha}$, $L_{\alpha}$, and $H_{\alpha}$ are shown in \cite[Proposition 3.4]{average}.
\bigskip

\noindent We present now the definitions of the one-stage optimization operators.
\begin{definition} Let $\alpha\in \RR_{+}$, $\rho\in \RR$, and $h\in \mathbb{M}(E)$.
Assume that for any $x\in E$ and $\Upsilon\in \mathbb{V}(x)$, $-\rho \mathcal{L}_{\alpha}(x,\Upsilon) +L_{\alpha}f(x,\Upsilon)+H_{\alpha}r(x,\Upsilon)+G_{\alpha} h(x,\Upsilon) $ is well defined.
The (ordinary) one-stage optimization operator is defined by
\begin{equation*}
\mathcal{T}_{\alpha}(\rho,h)(x)  =  \inf_{\Upsilon\in \mathbb{V}(x)} \Bigl\{-\rho \mathcal{L}_{\alpha}(x,\Upsilon)
+L_{\alpha}f(x,\Upsilon)+H_{\alpha}r(x,\Upsilon)+G_{\alpha} h(x,\Upsilon)  \Bigr\}.
\end{equation*}
Assume that for any $x\in E$ and $\Theta\in \mathbb{V}^{r}(x)$, $-\rho \mathcal{L}_{\alpha}(x,\Theta) +L_{\alpha}f(x,\Theta)+H_{\alpha}r(x,\Theta)+G_{\alpha} h(x,\Theta)$ is well defined.
The relaxed one-stage optimization operator is defined by
\begin{equation*}
\mathcal{R}_{\alpha}(\rho,h)(x)  = \inf_{\Theta\in \mathbb{V}^{r}(x)} \Bigl\{-\rho \mathcal{L}_{\alpha}(x,\Theta)
+L_{\alpha}f(x,\Theta)+H_{\alpha}r(x,\Theta)+G_{\alpha} h(x,\Theta) \Bigr\}.
\end{equation*}
\end{definition}

\noindent In particular for $\alpha =0$ we write for simplicity $\mathcal{T}_{0}=\mathcal{T}$, and $\mathcal{R}_{0}=\mathcal{R}$.

\bigskip

\noindent The sets of measurable selectors associated to $\bigl(\mathbb{U}(x)\bigr)_{x\in E}$, $\bigl(\mathbb{V}(x)\bigr)_{x\in E}$, $\bigl(\mathbb{V}^{r}(x)\bigr)_{x\in E}$ are defined by
$\mathcal{S}_{\mathbb{U}}  =  \Bigl\{  u \in \mathbb{M}(\widebar{E}, \mathbb{U}) : (\forall x\in \widebar{E}), u(x) \in \mathbb{U}(x)\Bigr\}$,
$\mathcal{S}_{\mathbb{V}}  =  \Bigl\{ (\nu,\nu_{\partial})\in \mathbb{M}(E, \mathbb{V}) : (\forall x\in E), \bigl( \nu(x), \: \nu_{\partial}(x) \bigr)\in \mathbb{V}(x)\Bigr\}$,
$\mathcal{S}_{\mathbb{V}^{r}} = \Bigl\{ (\mu,\mu_{\partial})\in \mathbb{M}(E, \mathbb{V}^{r}) : (\forall x\in E), \bigl( \mu(x), \: \mu_{\partial}(x) \bigr)\in \mathbb{V}^{r}(x)\Bigr\}$.

\bigskip

\noindent
\noindent For $\alpha\in \RR_{+}$, $\rho\in \RR$, and $v\in \mathbb{M}(E)$, the one-stage optimization problem associated to the operator $\mathcal{T}_{\alpha}(\rho,v)$, respectively
$\mathcal{R}_{\alpha}(\rho,v)$, consists of finding a measurable selector $\Upsilon\in \mathcal{S}_{\mathbb{V}}$, respectively $\Theta\in \mathcal{S}_{\mathbb{V}^{r}}$ such that for all $x\in E$,
$\mathcal{T}_{\alpha}(\rho,v)(x)  =  -\rho \mathcal{L}_{\alpha}(x,\Upsilon) + L_{\alpha}f(x,\Upsilon)+H_{\alpha}r(x,\Upsilon)+G_{\alpha} v(x,\Upsilon)$ and respectively
$\mathcal{R}_{\alpha}(\rho,v)(x)  =  -\rho \mathcal{L}_{\alpha}(x,\Theta) + L_{\alpha}f(x,\Theta)+H_{\alpha}r(x,\Theta)+G_{\alpha} v(x,\Theta)$.

\bigskip

\noindent Finally we conclude this section by recalling (see Propositions 3.8 and 3.10 in \cite{average}) that there exist two natural mappings from $\mathcal{S}_{\mathbb{U}}$ to $\mathcal{S}_{\mathbb{V}}$ and from $\mathcal{S}_{\mathbb{U}}$ to $\mathcal{U}$.
\begin{definition}
\label{mapu}
For $u \in\mathcal{S}_{\mathbb{U}}$, define the measurable mapping $u_{\phi}$ of the space $E$ into $\mathbb{V}$ by
\nl
$u_{\phi}$ $:$ $x$ $\rightarrow$ $\bigl(u(\phi(x,.)),u(\phi(x,t_{*}(x)))\bigr)$.
\end{definition}

\begin{definition}
\label{mapU}
For $u \in\mathcal{S}_{\mathbb{U}}$, define the measurable mapping $U_{u_{\phi}}$ of the space $\NN\times E\times \RR_{+}$ into $\mathbb{U}\times \mathbb{U}$
by $U_{u_{\phi}}$ $:$ $(n,x,t)$ $\rightarrow$ $\bigl(u(\phi(x,t)),u(\phi(x,t_{*}(x)))\bigr)$ of the space $\NN\times E\times \RR_{+}$ into $\mathbb{U}\times \mathbb{U}$.
\end{definition}

\begin{remark}
\label{feedbak} The measurable selectors of the kind $u_{\phi}$ as in Definition \ref{mapu} are called ordinary feedback measurable selectors in the class $\mathcal{S}_{\mathbb{V}}\subset \mathcal{S}_{\mathbb{V}^{r}}$
and the control strategies of the kind $U_{u_{\phi}}$ as in definition \ref{mapU} are called ordinary feedback control strategies in the class $\mathcal{U}$.
\end{remark}

\section{Assumptions}
\label{AssDef}
In order to prove our main results presented in section \ref{PIA}, we need to impose some conditions.
Assumptions \ref{Hyp3a}, \ref{Hyp6bis} and \ref{Hyp5a} are needed to guarantee some convergence and continuity properties of the one-stage optimization operators, and the existence of a measurable selector.
These properties are important to ensure the convergence of the policy iteration algorithm as shown in section \ref{ConvPIA}.
\begin{assumption}
\label{Hyp3a} For each $x\in E$, the restriction of $\lambda(x,.)$ to $\mathbb{U}(x)$ is continuous, for $t\in [0,t_{*}(x))$,
$\ds \int_{0}^{t} \sup_{a\in \mathbb{U}(\phi(x,s))} \lambda(\phi(x,s),a) \: ds < \infty$ and if $t_{*}(x)< \infty$ then $\ds \int_{0}^{t_{*}(x)} \sup_{a\in \mathbb{U}(\phi(x,s))} \lambda(\phi(x,s),a) \: ds < \infty$.
\end{assumption}

\begin{assumption}
\label{Hyp6bis} For all $y\in \widebar{E}$, the restriction of $f(y,.)$ to $ \mathbb{U}(y)$ is continuous and for all $z\in \partial E$,
the restriction of $r(z,.)$ to $\mathbb{U}(z)$ is continuous.
\end{assumption}

\begin{assumption}
\label{Hyp5a} For all $x \in \widebar{E}$ and $h\in \mathbb{B}(E)$, the restriction of $Qh(x,.)$ to $\mathbb{U}(x)$ is continuous.
\end{assumption}

\noindent
The next assumption is mainly used to show that the policy iteration algorithm converges to the optimal cost and gives an optimal feedback control as shown in section \ref{OptiPIA}.
This condition is somehow related to the so-called expected growth condition (see, for instance, Assumption 3.1 in \cite{guo06} for the discrete-time case, or
Assumption A in \cite{guo06a} for the continuous-time case).
\begin{assumption}
\label{A1}
Suppose that there exist $b\geq 0$, $c> 0$, $\delta>0$, $M\geq 0$ and $g\in\mathbb{M}^{ac}(E)$, $g\geq 1$
$\overline{r} \in\mathbb{M}(\partial E)$, $\overline{r}(z)\geq 0$,
satisfying for all $x\in E$
\begin{eqnarray}
& \ds \sup_{a\in \mathbb{U}(x)}\Bigl\{\mathcal{X}g(x)+c g(x)-\lambda(x,a)\left[g(x)-Qg(x,a)\right]\Bigr\} \leq b, &
\label{Cu1} \\
& \ds \sup_{a\in \mathbb{U}(x)}\Bigl\{f(x,a)\Bigr\}\leq M g(x), &
\label{Cu3}
\end{eqnarray}
and for all $x \in E$ with $t_{*}(x)<\infty$
\begin{eqnarray}
& \ds \sup_{a\in \mathbb{U}(\phi(x,t_{*}(x)))}\{\overline{r}(\phi(x,t_{*}(x)))+Qg(\phi(x,t_{*}(x)),a)\} \leq g(\phi(x,t_{*}(x))), &
\label{Cu2} \\
& \ds \sup_{a\in \mathbb{U}(\phi(x,t_{*}(x)))}\Bigl\{r(\phi(x,t_{*}(x)),a)\Bigr\}\leq \frac{M}{c+\delta}\overline{r} (\phi(x,t_{*}(x))).&
\label{Cu3a}
\end{eqnarray}
\end{assumption}

\noindent
In the next assumption notice that for any $u\in \mathcal{S}_\mathbb{U}$, $G(x,u_{\phi};.)$ can be seen as the stochastic kernel associated to the post-jump location of a PDMP.
This assumption is related to some geometric ergodic properties of the operator $G$ (see for example the comments on page 122 in \cite{hernandez99} or Lemma 3.3 in \cite{guo06}
for more details on this kind of assumption).
\begin{assumption}
\label{A2}
There exist $a> 0$, $0<\kappa<1$ and for any $u\in \mathcal{S}_\mathbb{U}$ there exists a probability measure $\nu_{u}$, such that $\nu_{u}(g)<+\infty$ and
\begin{eqnarray}
\bigl| G^kh(x,u_{\phi}) - \nu_{u}(h) \bigr| \leq a \|h\|_{g} \kappa^k g(x),
\label{Hy3}
\end{eqnarray}
for all $h\in \mathbb{B}_{g}(E)$ and $k\in\NN$.
\end{assumption}

\bigskip

\noindent
The following hypothesis is given by a Lyapunov-like inequality yielding an expected growth condition on the function $g$ with respect ot $G$ (for further comments on this kind of assumption, see for example section 10.2 in \cite[page 121]{hernandez99}).
\begin{assumption}
\label{A3}
There exist $0<k_{g}<1$ and $K_{g}\geq 0$ such that for all $x\in E$, $\Gamma \in \mathbb{V}(x)$,
\begin{eqnarray}
Gg(x,\Gamma) \leq k_{g} g(x) + K_{g}.
\label{eqA3}
\end{eqnarray}
\end{assumption}

\bigskip

\noindent
The final assumption is:
\begin{assumption}
\label{Hyp8a} There exist $\underline{\lambda} \in\mathbb{M}(E)^{+}$, and $K_{\lambda}\in \RR_{+}$ such that
\begin{enumerate}
\item [a)] $\lambda(y,a) \geq \underline{\lambda}(y)$ for all $y\in E$ and $a\in\mathbb{U}(y)$,
\item [b)] $\ds \int_0^{t_{*}(x)} e^{ct-\int_0^t \underline{\lambda}(\phi(x,s))ds}dt \leq K_{\lambda}$, for all $x\in E$,
\item [c)] $\ds \lim_{t\rightarrow +\infty} e^{ct-\int_0^{t}\underline{\lambda}(\phi(x,s))ds} =0$, for all $x\in E$ with $t_{*}(x)=+\infty$,
\item [d)] $\ds \lim_{t\rightarrow +\infty} e^{-\int_0^t \underline{\lambda}(\phi(x,s))ds} g(\phi(x,t)) =0$, for all $x\in E$ with $t_{*}(x)=\infty$,
\item [e)] $\ds \int_0^{t_{*}(x)} e^{-\int_0^t \underline{\lambda}(\phi(x,s))ds} \sup_{a\in \mathbb{U}(\phi(x,t))}f(\phi(x,t),a)dt < \infty$.
\end{enumerate}
\end{assumption}

\begin{remark} \label{vac} Notice the following consequences of Assumption \ref{Hyp8a}:
\begin{enumerate}
\item [i)] Assumption \ref{Hyp8a} c) implies that $\ds G_{\alpha}(x,\Theta;A) = \int_0^{t_{*}(x)}e^{-\alpha s - \Lambda^{\mu}(x,s)}\lambda QI_{A}(\phi(x,s),\mu(s)) ds$, and $H_{\alpha}w(x,\Theta) =0$,
for any $x\in E$ with $t_{*}(x)=+\infty$, $A\in \mathcal{B}(E)$, $\alpha\geq -c$, $\Theta=(\mu,\mu_{\partial})\in \mathbb{V}^{r}(x)$, $w\in \mathbb{M}(\partial E\times \mathbb{U})$.
\item [ii)] Assumptions \ref{Hyp8a} a) and b) imply that $\ds \mathcal{L}_{\alpha}(x,\Theta) \leq  K_{\lambda}$ for any $\alpha\geq -c$, $x\in E$, $\Theta\in \mathbb{V}^{r}(x)$.
\end{enumerate}
\end{remark}

\section{A pseudo-Poisson equation}
\label{pseudo}
We introduce in Definition \ref{psedoPoisson} a pseudo-Poisson equation associated to the stochastic kernel $G$.
Proposition \ref{PoissonRe} shows that there exists a solution for such an equation.
Moreover, it is proved in Proposition \ref{PoissonRe2} that this equation has the important characteristic of ensuring the policy improvement property
in the set $\mathcal{S}_\mathbb{U}$.

\begin{definition}
\label{psedoPoisson}
Consider $u\in \mathcal{S}_\mathbb{U}$. A pair $(\rho,h) \in \RR\times\mathbb{B}_{g}(E)$ is said to satisfies the pseudo-Poisson equation associated to $u$ if
\begin{align}
h(x) & = - \rho \mathcal{L}(x,u_{\phi}(x)) + Lf(x,u_{\phi}(x)) + Hr(x,u_{\phi}(x)) + Gh(x,u_{\phi}(x)).
\label{poisson}
\end{align}
\end{definition}

\begin{remark}
\label{difference}
This equation is clearly different from a classical Poisson equation encountered in the literature of the discrete-time Markov control processes
see for example equation (2.13) in \cite{PIA-97}.
In particular, the constant $\rho$, that will be shown to be the optimal cost, appears here as a multiplicative factor of the mapping $\mathcal{L}(x,u_{\phi}(x))$
and the costs $f$ and $r$ appear through the terms $Lf(x,u_{\phi}(x))$ and $Hr(x,u_{\phi}(x))$.
However, it will be shown in the following propositions that this pseudo-Poisson equation has still the good properties that we might expect to satisfy in order to guarantee
the convergence of the policy iteration algorithm.
\end{remark}

\begin{proposition}
\label{PoissonRe}
For arbitrary $u\in \mathcal{S}_\mathbb{U}$ the following assertions hold:
\begin{enumerate}
\item [(a)]
Set $\ds D_u=\int_{E} \mathcal{L}(y,u_{\phi}(y))\nu_{u}(dy)$. Then $0<D_u\leq K_{\lambda}$.
\item [(b)]
If $v\in \mathbb{B}_{g}(E)$ and $b\in \RR$ are such that for all $x\in E$,
\begin{equation}
v(x) = b \mathcal{L}(x,u_{\phi}(x)) + Gv(x,u_{\phi}(x)) \label{harm}
\end{equation}
then $b=0$ and for some $c_0\in \RR$, $v(x)=c_0$ for all $x\in E$.
\item [(c)]
Let $w_{u}$ be the mapping in $\mathbb{M}(E)$ defined by $w_{u}(x)  = Lf(x,u_{\phi}(x)) + Hr(x,u_{\phi}(x)) - \rho_u \mathcal{L}(x,u_{\phi}(x))$ for $x\in E$.
Define $(\rho_u,h_u)$ by
\begin{align}
\rho_u & = \frac{\ds \int_{E} \big[ Lf(y,u_{\phi}(y)) + Hr(y,u_{\phi}(y)) \big] \nu_{u}(dy)}{D_u} \geq 0,\label{rho1}\\
h_u(x) & = \sum_{k=0}^\infty G^kw_{u}(x,u_{\phi}(x)).
\label{hu}
\end{align}
Then $(\rho_u,h_u) \in \RR\times\mathbb{B}_{g}(E)$ and it is the unique solution to the Poisson equation (\ref{poisson}) associated to $u$
that satisfies
\begin{align}
\nu_{u}(h_u) = 0. \label{nuh0}
\end{align}
Moreover
\begin{align}
&\| h_u\|_g \leq  \frac{a M_u}{1-\kappa}, \text{ with } M_u := \max \Big\{ \rho_u K_{\lambda},\frac{M(1+bK_{\lambda})}{c} \Big\}.
\label{normh}
\end{align}
\end{enumerate}
\end{proposition}
\noindent \textbf{Proof:}
Item (a) is straightforward since $0<\mathcal{L}(x,u_{\phi}(x))\leq K_{\lambda}$ for all $x\in E$ (see Remark \ref{vac} ii)).

\bigskip

\noindent
For (b) let us suppose that $b\geq 0$. Since $0<\mathcal{L}(x,u_{\phi}(x))$ for all $x\in E$ it follows from (\ref{harm}) that
$v(x) \geq Gv(x,u_{\phi}(x))$ for all $x\in E$ and from Lemma 4.1 (a) in \cite{PIA-97}, $v(x)=c_0$ $\nu_{u}$-a.s. for some $c_0\in \RR$.
Returning to (\ref{harm}) and integrating with respect to $\nu_{u}$ we have that $0= b D_{u}$ and so $b=0$.
Therefore from (\ref{harm}), $v(x) = Gv(x,u_{\phi}(x))$, that is, $v$ is an $\nu_{u}$-harmonic function and therefore $v(x)=c_0$ for all $x\in E$
(see Lemma 4.1 (a) in \cite{PIA-97}).
If $b<0$ then from (\ref{harm}) it follows that $v(x) \leq Gv(x,u_{\phi}(x))$ for all $x\in E$ and from Lemma 4.1 (a) in \cite{PIA-97}, $v(x)=c_0$ $\nu_{u}$-a.s. for some $c_0\in \RR$.
Returning to (\ref{harm}) and integrating with respect to $\nu_{u}$
we have that $0= b D_{u}$ and since $D_{u}>0$, we have a contradiction.

\bigskip

\noindent For (c) we first note that from Proposition 3.12 in \cite{vanishing},
$0\leq Lf(x,u_{\phi}(x)) + Hr(x,u_{\phi}(x)) \leq \frac{M(1+bK_{\lambda})}{c}g(x)$ so that clearly
$\ds \int_{E} \big[ Lf(y,u_{\phi}(y)) + Hr(y,u_{\phi}(y)) \big] \nu_{u}(dy) < +\infty$, and thus (\ref{rho1}) is well defined.
Moreover $0\leq \rho_u \mathcal{L}(x,u_{\phi}(x))\leq \rho_u K_{\lambda}$ and thus $w_{u} \in \mathbb{B}_{g}(E)$ with
$\|w_{u}\|_g \leq M_u$ where $M_u$ is defined in (\ref{normh}). We also have from (\ref{rho1}) that
\begin{align}
\int_{E} w_{u}(y) \nu_{u}(dy) & =  \int_{E} \big[ Lf(y,u_{\phi}(y)) + Hr(y,u_{\phi}(y)) \big] \nu_{u}(dy) - \rho_u D_{u} \nonumber\\
&=  0
\label{zero}
\end{align}
and thus, from (\ref{Hy3}),
\begin{align}
\bigl| G^kw_{u}(x,u_{\phi}(x))\bigr| & = \bigl| G^kw_{u}(x,u_{\phi}(x)) - \nu_{u}(w_{u})\bigr| \leq a M_u  \kappa^k g(x),
\label{eqh1}
\end{align}
for all $x\in E$ and $k\in\NN$. From (\ref{hu}) and (\ref{eqh1}) it is clear that
\begin{equation}
\bigl| h_u(x) \bigr| \leq  \frac{a M_u}{1-\kappa} g(x), \label{eqh2}
\end{equation}
showing that $h_u$ is in $\mathbb{B}_{g}(E)$ and satisfies (\ref{normh}).
We also have from (\ref{hu}) that
\begin{align*}
h_u(x) - w_{u}(x) = \sum_{k=1}^\infty G^kw_{u}(x,u_{\phi}(x))= G_u h_u (x,u_{\phi}(x))
\end{align*}
showing that $(\rho_u,h_u)\in \RR\times\mathbb{B}_{g}(E)$ satisfies (\ref{poisson}).
\nl
If $(\rho_i,h_i)\in\RR\times\mathbb{B}_{g}(E)$, $i=1,2$
are 2 solutions to the Poisson equation (\ref{poisson}) then setting $v=h_1-h_2$ and $b=\rho_2-\rho_1$ we get that (\ref{harm}) is satisfied and
uniqueness follows from (b).
\hfill $\Box$

\bigskip

From now on, $(\rho_u,h_u)$ will denote the unique solution of the pseudo-Poisson equation (\ref{poisson}) that satisfies $\nu_{u}(h_{u}) = 0$.

\bigskip

The properties given in the following proposition are important for showing the convergence of the PIA.
\begin{proposition}
\label{PoissonRe2}
Consider $u\in \mathcal{S}_\mathbb{U}$.
Then there exists $\widehat{u}\in \mathcal{S}_\mathbb{U}$ such that
\begin{eqnarray}
\mathcal{R}(\rho_u,h_u)(x)  & = & -\rho_u \mathcal{L}(x,\widehat{u}_\phi(x))+Lf(x,\widehat{u}_\phi(x))+Hr(x,\widehat{u}_\phi(x))+G h_u(x,\widehat{u}_\phi(x)),
\label{eqPIA1}
\end{eqnarray}
and $\rho_{\widehat{u}}  \leq  \rho_u$.
\end{proposition}
\noindent \textbf{Proof:} From Theorem 3.22 in \cite{vanishing} we have that there exists $\widehat{u}\in \mathcal{S}_\mathbb{U}$ such that (\ref{eqPIA1})
holds. Clearly we have for every $x\in E$ that $h_u(x)\geq \mathcal{R}(\rho_u,h_u)(x)$, that is, from (\ref{eqPIA1}),
\begin{equation*}
h_u(x) \geq -\rho_u \mathcal{L}(x,\widehat{u}_\phi(x))+Lf(x,\widehat{u}_\phi(x))+Hr(x,\widehat{u}_\phi(x))+G h_u(x,\widehat{u}_\phi(x)).
\end{equation*}
Integrating the previous equation with respect to $\nu_{\widehat{u}}$ and recalling that the definition of $D_{u}$ (see item $a)$ in Proposition \ref{PoissonRe})
and $\ds \int_{E} Gh_u(y,\widehat{u}_\phi(y)) \nu_{\widehat{u}}(dy) =  \int_{E} h_u(y) \nu_{\widehat{u}}(dy)$, we get that
\begin{equation*}
\int_{E} h_u(y) \nu_{\widehat{u}}(dy) \geq  - \rho_u D_{\widehat{u}} + \rho_{\widehat{u}} D_{\widehat{u}}+\int_{E} h_u(y) \nu_{\widehat{u}}(dy)
\end{equation*}
that is, $\rho_u D_{\widehat{u}}\geq \rho_{\widehat{u}} D_{\widehat{u}}$ and since $D_{\widehat{u}}>0$ we get that $\rho_u \geq \rho_{\widehat{u}}$.
\hfill $\Box$

\section{The Policy Iteration Algorithm}
\label{PIA}
Having studied the pseudo-Poisson equation defined in section \ref{pseudo}, we are now in position to analyze the policy iteration algorithm.
In the first part, it is shown that the convergence of the policy iteration algorithm holds under a classical hypothesis (see for example assumption (H1) of Theorem 4.3 in \cite{PIA-97}).
Roughly speaking, it means that if the PIA computes a solution $(\rho_n,h_n)$ at the $n$th step then $(\rho_n,h_n)\rightarrow (\rho,h)$ and $(\rho,h)$ satisfies the optimality equation (\ref{OpEq}).
However it is far from obvious to claim that $\rho$ is actually the optimal cost for the long run average cost problem of the PDMP $\{X(t)\}$ and that there exists an optimal control.
In the second part of this section, these two issues are studied. In particular, we show that $\ds \rho= \inf_{U\in\mathcal{U}} \mathcal{A}(U,x)$ and
 the measurable selector $\widehat{u}_{\phi}$ of the optimality equation (\ref{OpEq}) provides an optimal control of the feedback form $U_{\widehat{u}_{\phi}}$ for the process $\{X(t)\}$:
$\ds \inf_{U\in\mathcal{U}} \mathcal{A}(U,x) = \mathcal{A}(U_{\widehat{u}_{\phi}},x)$.

\bigskip

\noindent
The policy iteration algorithm performs the following steps:
\begin{enumerate}
\item [Step 1:] Initialize with an arbitrary $u_0\in \mathcal{S}_\mathbb{U}$, and set $n=0$.
\item [Step 2:] Policy Evaluation - At the $n^{th}$-iteration consider $u_n\in \mathcal{S}_\mathbb{U}$
and evaluate $(\rho_n,h_n) \in \RR\times\mathbb{B}_{g}(E)$ the (unique) solution of the
Poisson equation (\ref{poisson}), (\ref{nuh0}) given by (\ref{rho1}) and (\ref{hu}), replacing $u$ by $u_n$,
thus we have that
\begin{align}
h_n(x) & = - \rho_n \mathcal{L}(x,(u_{n})_{\phi}(x)) + Lf(x,(u_{n})_{\phi}(x)) + Hr(x,(u_{n})_{\phi}(x)) + Gh_n(x,(u_{n})_{\phi}(x)), \label{poisson2}
\end{align}
with $\nu_{u_n}(h_n) = 0$.
\item [Step 3:] Policy Improvement - Determine $u_{n+1}\in \mathcal{S}_\mathbb{U}$ such that
\begin{eqnarray}
\mathcal{R}(\rho_n,h_n)(x) & = & -\rho_n \mathcal{L}(x,(u_{n+1})_{\phi}(x)) +Lf(x,(u_{n+1})_{\phi}(x))+Hr(x,(u_{n+1})_{\phi}(x)) \nonumber \\
& & +G h_n(x,(u_{n+1})_{\phi}(x)).
\label{eqPIA3}
\end{eqnarray}
\end{enumerate}

\bigskip

Notice that from Propositions \ref{PoissonRe} and \ref{PoissonRe2} the sequence $(\rho_n,h_n) \in \RR\times\mathbb{B}_{g}(E)$ and
$u_n\in \mathcal{S}_\mathbb{U}$ is well defined and moreover, $\rho_n\geq \rho_{n+1}\geq 0$. We set $\rho = \lim_{n\rightarrow \infty} \rho_n$.

\subsection{Convergence of the PIA}
\label{ConvPIA}
First we present in the next result some convergence properties of $G$, $H$, $L$ and $\mathcal{L}$.
\begin{proposition}
\label{lemlsc1}
Consider $h\in \mathbb{B}_{g}(E)$ and a sequence of functions $\big(h_{k}\big)_{k\in\NN}\in \mathbb{B}_{g}(E)$ such that for all $x\in E$,
$\ds \lim_{k\rightarrow \infty} h_{k}(x)=h(x)$ and
there exists $K_{h}$ satisfying $\bigl| h_{k}(x) \bigr| \leq K_{h} g(x)$ for all $k$ and all $x\in E$.
For $x\in E$, consider $\Theta_n=\bigl(\mu_n,\mu_{\partial,n}\bigr)\in \mathbb{V}^{r}(x)$ and $\Theta=\bigl(\mu,\mu_{\partial}\bigr)\in \mathbb{V}^{r}(x)$ such that
$\Theta_n \rightarrow \Theta$.
We have the following results:
\begin{align*}
&a) \lim_{n\rightarrow \infty}\mathcal{L}(x,\Theta_n) =\mathcal{L}(x,\Theta),&   &b) \lim_{n\rightarrow \infty}Lf(x,\Theta_n) = Lf(x,\Theta), \\
&c) \lim_{n\rightarrow \infty}Hr(x,\Theta_n) = Hr(x,\Theta), &   &d) \lim_{n\rightarrow \infty} Gh_{n}(x,\Theta_n) = Gh(x,\Theta).
\end{align*}
\end{proposition}
\textbf{Proof:}
The proof of item $a)$ is the same as in Proposition 5.7 in \cite{average} and it is essentially based on the fact that $\ds \lim_{n\rightarrow \infty} \Lambda^{\mu_{n}}(x,t)=\Lambda^{\mu}(x,t)$
by using assumption \ref{Hyp3a}.

\bigskip

\noindent \underline{Item b)} We have for $x\in E$,
\begin{eqnarray*}
L f(x,\Theta_n) & = & \int_0^{t_{*}(x)} \bigl[e^{-\Lambda^{\mu_{n}}(x,t)} - e^{-\Lambda^{\mu}(x,t)} \bigr] f(\phi(x,t),\mu_n(t))dt \\
&  & + \int_0^{t_{*}(x)} e^{-\Lambda^{\mu}(x,t)}f(\phi(x,t),\mu_n(s))dt.
\end{eqnarray*}
By combining items  $a)$ and $e)$ of assumption \ref{Hyp8a} and the dominated convergence theorem we obtain
 \begin{eqnarray*}
\lim_{n\rightarrow \infty}\int_0^{t_{*}(x)} \bigl| e^{-\Lambda^{\mu_{n}}(x,t)} - e^{-\Lambda^{\mu}(x,t)} \bigr| f(\phi(x,t),\mu_n(t))dt = 0.
\end{eqnarray*}
Therefore, we obtain item $b)$ by using assumption \ref{Hyp6bis}.

\bigskip

\noindent \underline{Item c)} Let us consider first that $t_{*}(x)= \infty$. From item $i)$ of remark \ref{vac}
it follows that $Hr(x,\Theta_n)= Hr(x,\Theta) = 0$.
Suppose now that $t_{*}(x) < \infty$ and set $z=\phi(x,t_{*}(x))$. From assumption \ref{Hyp6bis}, it follows that
$\ds \lim_{n\rightarrow \infty} r(z,\mu_{\partial,n}) = r(z,\mu_{\partial})$ showing item $c)$.

\bigskip

\noindent \underline{Item d)}
Let $\{\alpha_k\}$ a non increasing sequence of positive numbers with $\alpha_k\downarrow 0$. We have clearly
$\ds \liminf_{n\rightarrow \infty} Gh_{n}(x,\Theta_n) \geq \liminf_{n\rightarrow \infty}G_{\alpha_{n}}h_{n}(x,\Theta)$.
It follows that $\ds  \liminf_{n\rightarrow \infty} Gh_{n}(x,\Theta_n) \geq Gh(x,\Theta)$ by applying Proposition 3.18 in \cite{vanishing}.
Replacing $h_{n}$ by $-h_{n}$ it gives that $\ds \limsup_{n\rightarrow \infty} Gh_{n}(x,\Theta_n) \leq Gh(x,\Theta)$, completing the proof of item $d)$.
\hfill $\Box$

\bigskip

\noindent
We shall consider now the following assumption.
\begin{assumption}
\label{PIAConv} There exists a subsequence $\{h_k\}$ of $\{h_n\}$ and $h\in \mathbb{M}(E)$ such that for each $x\in E$,
\begin{equation}
\lim_{k\rightarrow  \infty} h_k(x)= h(x). \label{conver1}
\end{equation}
\end{assumption}

\bigskip

\noindent
The following theorem is the main result of this subsection. It shows the convergence of the PIA and ensures the existence of a measurable selector for the optimality equation.
\begin{theorem}
We have that $(\rho,h)\in \RR\times\mathbb{B}_{g}(E)$ satisfies the optimality equation:
\begin{equation}
h(x)= \mathcal{R}(\rho,h)(x).
\label{OpEq}
\end{equation}
Moreover there exists $\widehat{u}\in \mathcal{S}_\mathbb{U}$ such that
\begin{equation}
h(x)  = -\rho \mathcal{L}(x,\widehat{u}_\phi(x))+Lf(x,\widehat{u}_\phi(x))+Hr(x,\widehat{u}_\phi(x))+G h(x,\widehat{u}_\phi(x)).
\label{eqPIA1a}
\end{equation}
\end{theorem}
\noindent \textbf{Proof:} From (\ref{normh}) and recalling that $\rho_n\geq \rho_{n+1}$ we get that for all $k$,
\begin{align}
&\| h_k\|_g \leq \widetilde{M} := \frac{a M_{u_0}}{1-\kappa},\,\,\, M_{u_0} := \max\{\rho_0 K_{\lambda},\frac{M(1+bK_{\lambda})}{c}\}.\label{Mua1}
\end{align}
From (\ref{Mua1}) we get that $h\in \mathbb{B}_g(E)$, where $h$ is as in (\ref{conver1}). Consider $u_k\in \mathcal{S}_\mathbb{U}$ the
measurable selector associated to $(\rho_k,h_k)$ as in (\ref{poisson2}).
We have that for each $x\in E$,
$\mathbb{V}^{r}(x) $ is compact and $\{(u_k)_\phi\}$ is a sequence in $\mathcal{S}_{\mathbb{V}^{r}}$. Then according to Proposition 8.3 in \cite{PIA-97}
(see also \cite{schal75}) there exists $\Theta \in \mathcal{S}_{\mathbb{V}^{r}}$ such that $\Theta(x)\in \mathbb{V}^{r}(x)$ is an accumulation point
of $\{(u_k)_\phi(x)\}$ for each $x\in E$. Therefore for every $x\in E$, there exists a subsequence $k_i = k_i(x)$ such that $\lim_{i\rightarrow \infty}
(u_{k_i})_\phi(x) = \Theta(x)$. We fix now $x\in E$ and we consider the sub-sequence $k_i = k_i(x)$ as above. From Proposition \ref{lemlsc1} and
taking the limit in (\ref{poisson2}) for $n=k_i$ as $i\rightarrow \infty$  we have that
\begin{align}
h(x) & = - \rho \mathcal{L}(x,\Theta(x)) + Lf(x,\Theta(x)) + Hr(x,\Theta(x)) + Gh(x,\Theta(x)),\label{poisson4}
\end{align}
and thus clearly $h(x)\geq \mathcal{R}(\rho,h)(x)$. On the other hand from (\ref{poisson2}) and (\ref{eqPIA3}) we have that
\begin{align}
\mathcal{R}(\rho_{n-1},h_{n-1})(x) & +(\rho_{n-1}-\rho_{n}) \mathcal{L}(x,(u_{n})_{\phi}(x)) + G(h_n-h_{n-1})(x,(u_{n})_{\phi}(x))\nonumber\\
&= -\rho_{n} \mathcal{L}(x,(u_{n})_{\phi}(x))+Lf(x,(u_{n})_{\phi}(x))+Hr(x,(u_{n})_{\phi}(x))+G h_{n}(x,(u_{n})_{\phi}(x))\nonumber\\
&= h_n(x).\label{eqhn1}
\end{align}
From (\ref{eqhn1}) it is immediate that for any $\widetilde{\Theta} \in \mathcal{S}_{\mathbb{V}^{r}}$
\begin{align}
h_n(x)&\leq -\rho_{n-1} \mathcal{L}(x,\widetilde{\Theta}(x))+Lf(x,\widetilde{\Theta}(x))+Hr(x,\widetilde{\Theta}(x))+G h_{n-1}(x,\widetilde{\Theta}(x))\nonumber\\
& +(\rho_{n-1}-\rho_{n}) \mathcal{L}(x,(u_{n})_{\phi}(x)) + G(h_n-h_{n-1})(x,(u_{n})_{\phi}(x)).\label{eqhn2}
\end{align}
Fix $x$ and $k_i = k_i(x)$ as before and notice that for any $y\in E$, $\lim_{i \rightarrow \infty} (h_{k_i}(y)-h_{{k_i-1}}(y)) = 0$ and
from (\ref{Mua1}), $\| h_{k_i}-h_{{k_i-1}}\|_g \leq \widetilde{M}$. Applying Proposition \ref{lemlsc1} into (\ref{eqhn2}) replacing
$n$ by $k_i$ and taking the limit as $i\rightarrow \infty$ yields that
\begin{align}
h(x)&\leq -\rho \mathcal{L}(x,\widetilde{\Theta}(x))+Lf(x,\widetilde{\Theta}(x))+Hr(x,\widetilde{\Theta}(x))+Gh(x,\widetilde{\Theta}(x)),
\label{eqhn3}
\end{align}
and from (\ref{eqhn3}) we get that $h(x)\leq \mathcal{R}(\rho,h)(x)$. Thus we have (\ref{OpEq}).
\hfill $\Box$

\subsection{Optimality of the PIA}
\label{OptiPIA}
We present next a definition that will be useful for the next results.
\begin{definition}
\label{def1*} For any $\Theta=\bigl( \mu, \mu_{\partial} \bigr)\in \mathbb{V}$, define
\begin{eqnarray}
\bigl[\Theta\bigr]_{t}=\bigl( \mu(.+t), \mu_{\partial} \bigr).
\label{DefUps-rt}
\end{eqnarray}
\end{definition}

\bigskip

\noindent
Let us recall that the PDMP $\{\widehat{X}^{U}(t)\}$ and its associated components: $X(t)$, $Z(t)$, $N(t)$, $\tau(t)$  have been introduced in section \ref{pre} (see in particular equation (\ref{defXU})).
We need several auxiliary results (Propositions \ref{prop5bis}, \ref{prop5} and Corollary \ref{coro2b}) to show that the PIA actually provides an optimal solution for the average cost problem of the PDMP
$X(t)$.

\bigskip

\begin{proposition}
\label{prop5bis}
For $\hat{y}=(y,z,s,n)\in \widehat{E}$ and $U=(u,u_{\partial}) \in \mathbb{M}(\NN \times E\times \RR_{+};\mathbb{U}) \times \mathbb{M}(\NN \times E;\mathbb{U})$,
define $\Gamma^{U}(n,z)=\big(u(n,z,.),u_{\partial}(n,z)\big)\in \mathbb{V}$.
For $\epsilon \in (0,c)$ introduce
\begin{align}
\widehat{w}^{U}(\hat{y}) = & \widebar{c}L_{-\epsilon}f(y,\big[\Gamma^{U}(n,z)\big]_{s})+H_{-\epsilon}\widebar{r}(y,\big[\Gamma^{U}(n,z)\big]_{s})
+G_{-\epsilon}g(y,\big[\Gamma^{U}(n,z)\big]_{s}) \nonumber \\
& - b \mathcal{L}_{-\epsilon}(y,\big[\Gamma^{U}(n,z)\big]_{s}),
\end{align}
where $\widebar{c}= c - \epsilon$.
Then for all $x\in E$, $U\in \mathcal{U}$, we have
\begin{align}
E^{U}_{(x,0)} \Bigl[  \widehat{w}^{U}\bigl(\widehat{X}^{U}(t) \bigr)\Bigr] \leq e^{-\epsilon t}g(x)+\frac{b}{\epsilon}\bigl[1-e^{-\epsilon  t}\bigr].
\label{eqhgbis}
\end{align}
\end{proposition}
\noindent \textbf{Proof:}
For $\hat{y}=(y,z,s,n)\in \widehat{E}$ and $U=(u,u_{\partial}) \in \mathbb{M}(\NN \times E\times \RR_{+};\mathbb{U}) \times \mathbb{M}(\NN \times E;\mathbb{U})$,
define  $\widehat{f}^{U}(\hat{y})=f(y,u(n,z,s))$, $\widehat{r}^{U}(\hat{y})=\widebar{r}(y,u_{\partial}(n,z))$, $\widehat{g}(\hat{y})=g(y)$,
and for $t\in \RR_{+}$ $\widehat{\Lambda}^{U}(y,t)=\Lambda^{U}(x,n,t)$.
\nl
It is easy to show that $\widehat{w}^{U}\in \mathbb{M}(\widehat{E})$. Moreover, for $\hat{y}=(y,z,s,n)\in \widehat{E}$ and
$U=(u,u_{\partial}) \in \mathbb{M}(\NN \times E\times \RR_{+};\mathbb{U}) \times \mathbb{M}(\NN \times E;\mathbb{U})$, satisfying $\big[\Gamma^{U}(n,z)\big]_{s} \in \mathbb{V}(y)$
we have by using Corollary 3.11 in \cite{vanishing} with $\alpha=-\epsilon$ that
\begin{align}
\widebar{c}L_{-\epsilon}f(y,\big[\Gamma^{U}(n,z)\big]_{s})+H_{-\epsilon}\widebar{r}(y,\big[\Gamma^{U}(n,z)\big]_{s})
+G_{-\epsilon}g(y, & \big[\Gamma^{U}(n,z)\big]_{s}) \nonumber \\
&  - b \mathcal{L}_{-\epsilon}(y,\big[\Gamma^{U}(n,z)\big]_{s}) \leq g(y).
\label{ineqg}
\end{align}
Moreover, from Remark \ref{vac} $ii)$,
\begin{equation}
0< \mathcal{L}_{-\epsilon}(y,\big[\Gamma^{U}(n,z)\big]_{s}) \leq \mathcal{L}_{-c}(y,\big[\Gamma^{U}(n,z)\big]_{s}) \leq K_{\lambda}.
\label{Kxi}
\end{equation}
From now on, consider $U= (u,u_{\partial})\in \mathcal{U}$. Notice that for any $\hat{x}=(x,x,0,k)\in \widehat{E}$
\begin{eqnarray}
\widehat{w}^{U}(\widehat{x}) & = & \widebar{c}L_{-\epsilon}f(x,\Gamma^{U}(k,x))+H_{-\epsilon}\widebar{r}(y,\Gamma^{U}(k,x))
+G_{-\epsilon}g(x,\Gamma^{U}(k,x))  - b \mathcal{L}_{-\epsilon}(x,\Gamma^{U}(k,x)) \nonumber \\
& = & \int_{0}^{t_{*}(x)} e^{\epsilon s -\Lambda^{\nu_{k}}(x,s)}\biggl[ -b + \widebar{c} f(\phi(x,s),\nu_{k}(s)) + \lambda(\phi(x,s),\nu_{k}(s)) Qg(\phi(x,s),\nu_{k}(s))\biggr] ds \nonumber \\
& &+ e^{\epsilon t_{*}(x) -\Lambda^{\nu_{k}}(x,t_{*}(x))} \Bigl[ Qg(\phi(x,t_{*}(x)),u_{\partial}(k,x)) + \widebar{r}(\phi(x,t_{*}(x)),u_{\partial}(k,x)) \Bigr],
\label{ineqgb1}
\end{eqnarray}
with $\nu_{k}(.)=u(k,x,.)$.
Since for all $k\in \NN$, $x\in E$, $\Gamma^{U}(k,x)\in \mathbb{V}(x)$, it follows from equation (\ref{ineqg}) that
\begin{eqnarray}
\widehat{w}^{U}(\widehat{x}) & \leq & g(x).
\label{ineqgb}
\end{eqnarray}
Moreover, since $\big[\Gamma^{U}(N(t),Z(t))\big]_{\tau(t)} \in \mathbb{V}(X(t))$,
the inequality (\ref{Kxi}) implies that
\begin{align*}
 J^{U}_{m}(t,\hat{x}) &:=  E^{U}_{(x,k)} \Biggl[  \int_{0}^{t\wedge T_{m}} e^{\epsilon s} \Bigl[ \widebar{c} \widehat{f}^{U}(\widehat{X}^{U}(s)) -b \Bigr] ds
+ \int_{0}^{t\wedge T_{m}} e^{\epsilon s} \widehat{r}^{U}\bigl(\widehat{X}^{U}(s-) \bigr) dp^{*}(s) \nonumber \\
 &  + e^{\epsilon (t\wedge T_{m})} \widehat{w}^{U}\bigl(\widehat{X}^{U}(t\wedge T_{m}) \bigr) \biggr) \Biggr],
\end{align*}
is well defined for any $\hat{x}=(x,x,0,k)\in \widehat{E}$.
\nl
Let us show by induction on $m\in \NN$ that $J^{U}_{m}(t,\hat{x}) \leq g(x)$ for all $t\in \RR_{+}$, $\hat{x}=(x,x,0,k)\in \widehat{E}$.
Clearly, we have that $J^{U}_{0}(t,\hat{x})=\widehat{w}^{U}(\widehat{x})$.
Consequently, from equation (\ref{ineqgb}), we have that $J^{U}_{0}(t,\hat{x})\leq g(x)$
for all $t\in \RR_{+}$, $\hat{x}=(x,x,0,k)\in \widehat{E}$.
Now assume that for $m\in \NN$ we have that $J^{U}_{m}(t,\hat{x}) \leq g(x)$ for all $t\in \RR_{+}$, $\hat{x}=(x,x,0,k)\in \widehat{E}$.
Following the same arguments as in the proof of Proposition 4.3 in \cite{average}, it is easy to show that for $t\in \RR_{+}$
\begin{align}
 J^{U}_{m+1} & (t,\hat{x})  \leq \int_{0}^{t\wedge t_{*}(x)} e^{\epsilon s -\Lambda^{\nu_{k}}(x,s)}\biggl[ -b + \widebar{c} f(\phi(x,s),\nu_{k}(s)) + \lambda(\phi(x,s),\nu_{k}(s)) Qg(\phi(x,s),\nu_{k}(s))\biggr] ds \nonumber \\
& \phantom{=} + I_{\{t\geq t_{*(x)}\}} e^{\epsilon t_{*}(x) -\Lambda^{\nu_{k}}(x,t_{*}(x))} \Bigl[ Qg(\phi(x,t_{*}(x)),u_{\partial}(k,x)) +\widebar{r}(\phi(x,t_{*}(x)),u_{\partial}(k,x)) \Bigr] \nonumber \\
& \phantom{=}  +  I_{\{t<t_{*(x)}\}} e^{\epsilon t-\Lambda^{\nu_{k}}(x,t)}\widehat{w}^{U}(\widehat{\phi}(\hat{x},t)). \label{eq5prop2a}
\end{align}
Now if $t<t_{*}(x)$, then by using the fact that $\widehat{\phi}(\hat{x},t)=\big(\phi(x,t),x,t,k\big)$ we get that
\begin{eqnarray*}
\widehat{w}^{U}(\widehat{\phi}(\hat{x},t)) & = & \widebar{c}L_{-\epsilon}f(x,\big[\Gamma^{U}(k,x)\big]_{t})+H_{-\epsilon}\widebar{r}(x,\big[\Gamma^{U}(k,x)\big]_{t})
+G_{-\epsilon}g(x,\big[\Gamma^{U}(k,x)\big]_{t}) \nonumber \\
& & - b \mathcal{L}_{-\epsilon}(x,\big[\Gamma^{U}(k,x)\big]_{t}),
\end{eqnarray*}
and it follows, by applying Proposition 4.2 in \cite{average}, that
\begin{eqnarray}
\widehat{w}^{U}(\hat{x}) & = &\int_{0}^{t} e^{\epsilon s -\Lambda^{\nu_{k}}(x,s)}\biggl[ -b + \widebar{c} f(\phi(x,s),\nu_{k}(s))
+ \lambda(\phi(x,s),\nu_{k}(s)) Qg(\phi(x,s),\nu_{k}(s))\biggr] ds \nonumber \\
& & + e^{\epsilon t-\Lambda^{\nu_{k}}(x,t)}\widehat{w}^{U}(\widehat{\phi}(\hat{x},t)).
\label{eq5prop2ab}
\end{eqnarray}
Therefore, combining equations (\ref{eq5prop2a}) and (\ref{eq5prop2ab}) we get that $J^{U}_{m+1} (t,\hat{x}) \leq \widehat{w}^{U}(\hat{x})$ and by using equation (\ref{ineqgb}) we have that $J^{U}_{m+1} (t,\hat{x}) \leq g(x)$.
\nl
If $t\geq t_{*}(x)$, then equations (\ref{ineqgb1}) and (\ref{eq5prop2a}) yields $J^{U}_{m}(t,\hat{x}) \leq \widehat{w}^{U}(\hat{x})$.
By using equations (\ref{ineqgb}),  we have $J^{U}_{m}(t,\hat{x}) \leq g(x)$,
showing the fact that for all $m\in \NN$, $J^{U}_{m}(t,\hat{x}) \leq g(x)$ for all $t\in \RR_{+}$, $\hat{x}=(x,x,0,k)\in \widehat{E}$.
\nl
Consequently, this implies that
$\ds -b E^{U}_{(x,0)} \Bigl[ \int_{0}^{t\wedge T_{m}} e^{\epsilon s} ds\Bigr]
+ E^{U}_{(x,0)} \Bigl[ e^{\epsilon (t\wedge T_{m})} \widehat{w}^{U}\bigl(\widehat{X}^{U}(t\wedge T_{m}) \bigr)\Bigr] \leq g(x)$.
Combining Fatou's Lemma and Remark \ref{RemJump} we obtain that
\begin{align}
-\frac{b}{\epsilon}\bigl[e^{\epsilon t}-1 \bigr] + e^{\epsilon t} E^{U}_{(x,0)} \Bigl[  \widehat{w}^{U}\bigl(\widehat{X}^{U}(t) \bigr)\Bigr] \leq g(x),
\label{eqprop5}
\end{align}
showing the result.
 \hfill $\Box$

\bigskip

\noindent

\begin{proposition}
\label{prop5}
For all $x\in E$, $U\in \mathcal{U}$, we have that $\ds E^{U}_{(x,0)} \Bigl[  \widehat{w}^{U}\bigl(\widehat{X}^{U}(t\wedge T_{m}) \bigr)\Bigr]$ exists in $\RR_{+}$ for any $(t,m) \in \RR_{+}\times \NN$ and
\begin{eqnarray}
\limsup_{t\rightarrow +\infty} \frac{1}{t} \limsup_{m\rightarrow \infty} E^{U}_{(x,0)} \Bigl[  \widehat{w}^{U}\bigl(\widehat{X}^{U}(t\wedge T_{m}) \bigr)\Bigr] = 0.
\label{eqlim}
\end{eqnarray}
\end{proposition}
\noindent \textbf{Proof:}
Clearly, we have
\begin{eqnarray*}
E^{U}_{(x,0)} \Bigl[  \widehat{w}^{U}\bigl(\widehat{X}^{U}(t\wedge T_{m}) \bigr)\Bigr]
=E^{U}_{(x,0)} \Bigl[  I_{\{t<T_{m}\}}\widehat{w}^{U}\bigl(\widehat{X}^{U}(t) \bigr)\Bigr]
+ E^{U}_{(x,0)} \Bigl[  I_{\{t\geq T_{m}\}}\widehat{w}^{U}\bigl(\widehat{X}^{U}(T_{m}) \bigr)\Bigr] ,
\end{eqnarray*}
and thus by using Remark \ref{vac} $ii)$,
\begin{eqnarray}
0\leq E^{U}_{(x,0)} \Bigl[  \widehat{w}^{U} \bigl(\widehat{X}^{U}(t\wedge T_{m}) \bigr)\Bigr]
\leq E^{U}_{(x,0)} \Bigl[ \widehat{w}^{U}\bigl(\widehat{X}^{U}(t) \bigr)\Bigr]
+E^{U}_{(x,0)} \Bigl[ \widehat{w}^{U}\bigl(\widehat{X}^{U}(T_{m})\bigr)\Bigr] +bK_{\lambda}.
\label{eqhgbis2}
\end{eqnarray}
Iterating Assumption \ref{A3}, we obtain that for all $m\in \NN$,
$\ds E^{U}_{(x,0)} \Bigl[ \widehat{w}^{U}\bigl(\widehat{X}^{U}(T_{m})\bigr)\Bigr] \leq g(x) + \frac{K_{g}}{1-k_{g}}$.
Combining equations (\ref{eqhgbis}), (\ref{eqhgbis2}) and the previous inequality, the result follows.
\hfill$\Box$

\begin{corollary}
\label{coro2b}
For all $U\in \mathcal{U}$,
\begin{eqnarray}
\limsup_{t\rightarrow +\infty} \frac{1}{t} \limsup_{m\rightarrow \infty} E^{U}_{(x,0)} \Bigl[h\big(X(t\wedge T_{m})\big)\Bigr] \leq 0,
\label{eqlim2}
\end{eqnarray}
and
\begin{eqnarray}
\limsup_{t\rightarrow +\infty} \frac{1}{t} \limsup_{m\rightarrow \infty} E^{U_{\widehat{u}_{\phi}}}_{(x,0)} \Bigl[h\big(X(t\wedge T_{m})\big)\Bigr] = 0.
\label{eqlim3}
\end{eqnarray}
\end{corollary}
\noindent \textbf{Proof:}
From equation (\ref{eqPIA1a}), it follows that for all $x\in E$, $\Gamma \in \mathbb{V}(x)$,
\begin{eqnarray}
-\rho \mathcal{L}(x,\widehat{u}_{\phi}(x))+Gh(x,\widehat{u}_{\phi}(x)) \leq h(x) \leq Lf(x,\Gamma)+Hr(x,\Gamma)+Gh(x,\Gamma).
\label{coro1}
\end{eqnarray}
Consequently, by using Remark \ref{vac} $ii)$, the definition of $\widehat{w}$ and Assumption \ref{A1} we obtain that there exist $M_{1}>0$ such that for any $U\in \mathcal{U}$
\begin{eqnarray*}
h\big(X(t\wedge T_{m})\big)&\leq& M_{1} \Big[ \widehat{w}^{U} \bigl(\widehat{X}^{U}(t\wedge T_{m}) \bigr) + bK_{\lambda} \Big].
\end{eqnarray*}
Consequently, combining the previous equation and (\ref{eqlim}) we obtain equation (\ref{eqlim2}).
\nl
Moreover, notice that $\big[\Gamma^{U_{\widehat{u}_{\phi}}}(N(t),Z(t))\big]_{\tau(t)}=\widehat{u}_{\phi}(X(t))$ and so equation (\ref{coro1}) implies
\begin{eqnarray*}
-\|h\|_{g} \Big[ \widehat{w}^{U_{\widehat{u}_{\phi}}} \bigl(\widehat{X}^{U_{\widehat{u}_{\phi}}}(t\wedge T_{m}) \bigr) + bK_{\lambda} \Big] -\rho K_{\lambda} \leq h\big(X(t\wedge T_{m})\big).
\end{eqnarray*}
By using equation (\ref{eqlim}), this yields that $\ds \limsup_{t\rightarrow +\infty} \frac{1}{t} \limsup_{m\rightarrow \infty} E^{U_{\widehat{u}_{\phi}}}_{(x,0)} \Bigl[h\big(X(t\wedge T_{m})\big)\Bigr] \geq 0$.
Combining the previous inequality with (\ref{eqlim2}), the result follows.
\hfill$\Box$

\bigskip

\noindent
Finally, we can now present our second main result. It states that the measurable selector $\widehat{u}_{\phi}$ of the optimality equation (\ref{OpEq}) associated to $(\rho,h)$ gives an optimal feedback control $U_{\widehat{u}_{\phi}}$ for the process $\{X(t)\}$.
\begin{theorem}
\label{prop2b}
The control $U_{\widehat{u}_{\phi}}$ is an optimal strategy for the long-run average control problem:
\begin{eqnarray*}
\rho= \inf_{U\in\mathcal{U}} \mathcal{A}(U,x) = \mathcal{A}(U_{\widehat{u}_{\phi}},x),
\end{eqnarray*}
for all $x\in E$.
\end{theorem}
\noindent \textbf{Proof:}
From Proposition \ref{prop5} we have that $\ds E^{U}_{(x,0)} \Bigl[h\big(X(t\wedge T_{m})\big)\Bigr] = E^{U}_{(x,0)} \Bigl[h\big(X(t\wedge T_{m})\big)\Bigr]$ is well defined.
Therefore, following the same arguments as in Proposition 4.3 in \cite{average} it can be shown that
\begin{align*}
E^{U}_{(x,0)} \biggl[ & \int_{0}^{t\wedge T_{m}} f\bigl(X(s),u(N(s),Z(s),\tau(s)) \bigr) ds  + \int_{0}^{t\wedge T_{m}} r\bigl(X(s-),u_{\partial}(N(s),X(s-)) \bigr) dp^{*}(s)\biggr] \nonumber \\
& + E^{U}_{(x,0)} \Bigl[ h\bigl(X(t\wedge T_{m}) \bigr)\Bigr] \geq  E^{U}_{(x,0)} \Bigl[ \rho [t\wedge T_{m}] \Bigr] +h(x),
\end{align*}
where $U=\big(u,u_{\partial}\big)\in \mathcal{U}$.
From equation (\ref{eqlim2}), it implies that
\begin{align*}
\limsup_{t\rightarrow +\infty} \frac{1}{t} E^{U}_{(x,0)} \biggl[ & \int_{0}^{t} f\bigl(X(s),u(N(s),Z(s),\tau(s)) \bigr) ds  + \int_{0}^{t} r\bigl(X(s-),u_{\partial}(N(s),X(s-)) \bigr) dp^{*}(s)\biggr]  \geq  \rho,
\end{align*}
showing that $\ds \inf_{U\in\mathcal{U}} \mathcal{A}(U,x) \geq \rho$.
\nl
From equation (\ref{eqlim3}), it can be shown by using the same arguments as in the proof of Proposition 4.4 in \cite{average} that
\begin{align*}
\limsup_{t\rightarrow +\infty} \frac{1}{t} E^{U_{\widehat{u}_{\phi}}}_{(x,0)} \biggl[  \int_{0}^{t} f\bigl(X(s),\widehat{u}(X(s) &) \bigr) ds  + \int_{0}^{t} r\bigl(X(s-),\widehat{u}(X(s-)) \bigr) dp^{*}(s)\biggr] \nonumber \\
& \leq  \rho -\limsup_{t\rightarrow +\infty} \frac{1}{t} \limsup_{m\rightarrow \infty} E^{U_{\widehat{u}_{\phi}}}_{(x,0)} \Bigl[h\big(X(t\wedge T_{m})\big)\Bigr] = \rho,
\end{align*}
implying that $\ds \inf_{U\in\mathcal{U}} \mathcal{A}(U,x) \leq \rho$.
\nl
Therefore, it follows that $\ds \rho= \inf_{U\in\mathcal{U}} \mathcal{A}(U,x) = \mathcal{A}(U_{\widehat{u}_{\phi}},x)$ for all $x\in E$.
\hfill$\Box$

\bibliography{PIA}

\end{document}